\documentclass[a4paper,12pt,reqno]{amsart}
\usepackage{amssymb,amsthm}
\usepackage{ifthen}
\usepackage{graphicx}
\usepackage{float}
\usepackage{mathrsfs}

\theoremstyle{plain}

\newtheorem{thm}{Theorem}[section]
\newtheorem*{thm*}{Theorem}
\newtheorem{theorem}{Theorem}[section]
\numberwithin{equation}{section}

\newtheorem{cor}{Corollary}[section]
\newtheorem{lem}{Lemma}[section]
\newtheorem{prop}{Proposition}[section]
\newtheorem{rem}{Remark}[section]
\newtheorem{example}{Example}[section]
\newtheorem{defn}{Definition}[section]

\theoremstyle{definition}
\newcounter {own}
\def\theown {\thesection  .\arabic{own}}

\newenvironment{pf}[1][]{%
 \vskip 3mm
 \noindent
 \ifthenelse{\equal{#1}{}}%
  {{\slshape Proof. }}%
  {{\slshape #1.} }%
 }%
{\qed\bigskip}

\newcounter{alphabet}
\newcounter{tmp}

\newcommand{\ds}{\displaystyle}

\newcounter{minutes}
\divide\time by 60
\newcounter{hours}
\multiply\time by 60 \addtocounter{minutes}{-\time}

\begin{document}
\bibliographystyle{amsplain}
\title{Shift-invariant spaces with countably many mutually orthogonal generators on the Heisenberg group}

\thanks{
File:~\jobname .tex,
          printed: \number\year-\number\month-\number\day,
          \thehours.\ifnum\theminutes<10{0}\fi\theminutes}

\author{R. Radha $^\dagger$}

\address{R. Radha, Department of Mathematics,
Indian Institute of Technology Madras, Chennai--600 036, India.}
\email{radharam@iitm.ac.in}
\author{Saswata Adhikari }
\address{Saswata Adhikari, Department of Mathematics,
Indian Institute of Technology Madras, Chennai--600 036, India.}
\email{saswata.adhikari@gmail.com}
\keywords{Bessel sequence, frames, group Fourier transform, Heisenberg group, orthonormal system, twisted translation.\\
$^\dagger$ {\tt Corresponding author}
}

\maketitle
\pagestyle{myheadings}
\markboth{R. Radha and Saswata Adhikari}{Shift-invariant spaces with countably many mutually orthogonal generators on the Heisenberg group}
\begin{abstract}
Let $E(\mathscr{A})$ denote the shift-invariant space associated with a countable family $\mathscr{A}$ of functions in $L^{2}(\mathbb{H}^{n})$ with mutually orthogonal generators, where $\mathbb{H}^{n}$ denotes the Heisenberg group. The characterizations for the collection $E(\mathscr{A})$ to be orthonormal, Bessel sequence, Parseval frame and so on are obtained in terms of the group Fourier transform of the Heisenberg group. These results are derived using such type of results which were proved for twisted shift-invariant spaces and characterized in terms of Weyl transform. In the last section of the paper, some results on oblique dual of the left translates of a single function $\varphi$ is discussed in the context of principal shift-invariant space $V(\varphi)$.
\end{abstract}

\section{Introduction}
Let $\mathscr{L}$ be a discrete subgroup of the Heisenberg group $\mathbb{H}^{n}$ such that $\mathbb{H}^{n}/\mathscr{L}$ is compact. In other words, $\mathscr{L}$ is a lattice in $\mathbb{H}^{n}$. For $\varphi\in L^{2}(\mathbb{H}^{n})$, the principal shift-invariant space, denoted by $V(\varphi)$, is defined to be  $\overline{span}\{{L_{l}\varphi:l\in\mathscr{L}}\}$, where
$L_{l}\varphi(X)=\varphi(l^{-1}.X),~X\in\mathbb{H}^{n}$. Without loss of generality and to avoid complexity in proofs, we shall work with the standard lattice
$\mathscr{L}=\{(2k,l,m):k,l\in\mathbb{Z}^{n},m\in\mathbb{Z}\}$.

Let $\mathscr{A}=\{\varphi_{j}:j\in\mathbb{N}\}$ be a countable family of functions in $L^{2}(\mathbb{H}^{n})$ such that $V(\varphi_{i})$ is orthogonal to $V(\varphi_{j})$ whenever $i\neq j$. Let $E(\mathscr{A})$ denote the collection $\{L_{(2k,l,m)}\varphi_{j}: (k,l,m)\in\mathbb{Z}^{2n+1}, j\in\mathbb{N}\}$. We shall denote span of $E(\mathscr{A})$ by $U(\mathscr{A})$ and  $\overline{span}~E(\mathscr{A})$ by $V(\mathscr{A})$.

The aim of this paper is to obtain characterizations for the collection $E(\mathscr{A})$ to be an orthonormal system, Bessel sequence, Parseval frame and so on in terms of the group Fourier transform on the Heisenberg group. The proofs of these results are obtained by using such type of results which were proved for twisted shift-invariant spaces and characterized in terms of the kernel of the Weyl transform.  We  refer to \cite{ras} of the authors in this connection. As in \cite{bo}, we introduce a map $\mathscr{T}$ from $L^{2}(\mathbb{H}^{n})$ into $L^{2}((0,1],\ell^{2}(\mathbb{Z},\mathcal{B}_{2}))$, using the group Fourier transform on $\mathbb{H}^{n}$, which turns out to be an isometric isomorphism. Using this map $\mathscr{T}$, we prove that $E(\mathscr{A})$ is an orthonormal system in $L^{2}(\mathbb{H}^{n})$ iff $\{\mathscr{T}\varphi_{j}(\lambda):j\in\mathbb{N}\}$ is an orthonormal system in $\ell^{2}(\mathbb{Z},\mathcal{B}_{2}))$ for a.e. $\lambda\in (0,1]$ and each $\varphi_{j}$ satisfies an additional condition, which we call ``condition C". We also obtain similar type of characterizations for the collection $E(\mathscr{A})$ to be Bessel sequence, Parseval frame and so on.

In the last part of this paper, we study in detail on frames for the principal shift-invariant space $V(\varphi)$. Given two Bessel sequences $\{L_{(2k,l,m)}\varphi:(k,l,m)\in\mathbb{Z}^{2n+1}\}$ and $\{L_{(2k,l,m)}\tilde{\varphi}:(k,l,m)\in\mathbb{Z}^{2n+1}\}$, we provide a necessary and sufficient condition on the generators $\varphi$ and $\tilde{\varphi}$ such that $\{L_{(2k,l,m)}\tilde{\varphi}:(k,l,m)\in\mathbb{Z}^{2n+1}\}$ is an oblique dual of $\{L_{(2k,l,m)}\varphi:(k,l,m)\in\mathbb{Z}^{2n+1}\}$ under a certain additional assumption. We also show that if $\{L_{(2k,l,m)}\varphi:(k,l,m)\in\mathbb{Z}^{2n+1}\}$ is a frame for $V(\varphi)$, then there exists $\tilde{\varphi}\in V(\varphi)$ such that $\{L_{(2k,l,m)}\tilde{\varphi}:(k,l,m)\in\mathbb{Z}^{2n+1}\}$ is an oblique dual of
$\{L_{(2k,l,m)}\varphi:(k,l,m)\in\mathbb{Z}^{2n+1}\}$.

Characterizations of shift-invariant spaces in $L^{2}(\mathbb{R}^{n})$ in terms of the range function were studied by Bownik in \cite{bo}. The results in this paper were generalized to locally compact abelian group by Cabrelli and Paternostro in \cite{cab}. Shift-invariant spaces on locally compact abelian groups were also independently studied by Gol and Tousi in \cite{kgr}. Radha and Kumar studied shift-invariant spaces for a compact non-abelian group in \cite{rad}. The concept of the bracket map has been generalized in \cite{bhp} to include non-abelian discrete group $\Gamma$ using its unitary representations and $L^{1}$ space over the non commutative measurable space vNa$(\Gamma)$, which is the compact dual of $\Gamma$ whose underlying space is a group von Neumann algebra. Using this bracket map characterizations of orthonormal basis, Riesz basis, frames were obtained for shift-invariant spaces in a Hilbert space $\mathcal{H}$ given by the action of a non-abelian countable discrete group $\Gamma$. In \cite{lf}, Luef provided a connection between the construction of projections in non-commutative tori and the construction of tight Gabor frames for $L^{2}(\mathbb{R})$.

In \cite{cmo}, shift-invariant spaces are characterized in terms of the range function for $SI/Z$ type groups. These are simply connected Lie groups having unitary irreducible representations that are square-integrable modulo the center. Using the range function $\mathscr{T}$, they obtain characterizations for  $E(\mathscr{A})$ to be frame, Riesz basis in the following sense. They show that $E(\mathscr{A})$ is a frame iff $\{\mathscr{T}(L_{k}\varphi)(\sigma):\varphi\in\mathscr{A},k\in\Gamma_{1}\}$ is a frame for $J(\sigma)$ a.e. $\sigma\in\mathbb{T}^{n}$, where $\mathscr{T}$ is a range function. This problem is different from that of \cite{bo}, where Bownik considered the group $\mathbb{R}^{n}$ and showed that $E(\mathscr{A})$ is a frame iff $\{\mathscr{T}\varphi(x):\varphi\in\mathscr{A}\}$ is a frame for $J(x)$ a.e. $x\in\mathbb{T}^{n}$. In this paper, we also wish to consider a similar type of problem as in \cite{bo} in the case of Heisenberg group. We obtain the results with an additional assumption on $\mathscr{A}=\{\varphi_{j}:j\in\mathbb{N}\}$ that $V(\varphi_{i})$ is orthogonal to $V(\varphi_{j})$ whenever $i\neq j$. We want to emphasize that the technicality used in the proofs of our theorems are different from that of \cite{bo} as well as \cite{cmo}.

We also would like to mention here that, when $\mathscr{A}=\{\varphi\}$, our result for $E(\mathscr{A})$ to be an orthonormal system is similar to that of \cite{bhm} proved for the polarized Heisenberg group. But our approach and proofs are totally different from that of \cite{bhm}. We make use of explicit computations using the Weyl transform, twisted shift-invariant spaces and periodization associated with projective representation of $\mathbb{C}^{n}$ instead of looking into integer periodization in the unitary dual of the group which corresponds to the central variable as in \cite{bhm}.

The Heisenberg group $\mathbb{H}^n$ is a simple example of a non-abelian, non-compact group. It is a nilpotent Lie group whose underlying manifold is $\mathbb{C}^n\times\mathbb{R}$, where the group operation is defined by $(z,t).(w,s)=(z+w,t+s+\dfrac{1}{2}\mathrm{Im}z.\overline{w})$ and the Haar measure is the Lebesgue measure  $\mathrm{d}z\mathrm{d}t$ on $\mathbb{C}^n\times\mathbb{R}$. By Stone-von Neumann theorem, every infinite dimensional irreducible unitary representation on the Heisenberg group is unitarily equivalent to the representation $\pi_{\lambda}, \lambda\in \mathbb{R}^{\star}$, where $\pi_{\lambda}$ is defined by
\begin{eqnarray*}
\pi_{\lambda}(z,t)\varphi(\xi)=e^{2\pi i\lambda t}e^{2\pi i\lambda(x.\xi+\frac{1}{2}x.y)}\varphi(\xi+y),
\end{eqnarray*}
for $z=x+iy$ and $\varphi\in L^{2}(\mathbb{R}^{n})$.
The group Fourier transform on $\mathbb{H}^n$ is defined to be $$\widehat{f}(\lambda)=\int\limits_{\mathbb{H}^n}f(z,t)\pi_\lambda(z,t)dzdt,
\hspace{.5cm}\lambda\in\mathbb{R}^{\star},$$ for $f\in L^{1}(\mathbb{H}^n)$. Then $\widehat{f}(\lambda)\in\mathcal{B}(L^{2}(\mathbb{R}^{n}))$, the class of bounded operators on $L^{2}(\mathbb{R}^{n})$ and
\begin{eqnarray*}
\|\widehat{f}(\lambda)\|_{\mathcal{B}}\leq \|f\|_{L^{1}(\mathbb{H}^{n})}.
\end{eqnarray*}
Further, the group Fourier transform is an isometric isomorphism of $L^{2}(\mathbb{H}^{n})$ onto $L^{2}(\mathbb{R}^{\star},\mathcal{B}_{2};d\mu)$ , where $\mathcal{B}_{2}$ denotes the space of Hilbert-Schmidt operators on $L^{2}(\mathbb{R}^{n})$ and $d\mu(\lambda)=|\lambda|^{n}d\lambda$. More precisely, for $f\in L^{2}(\mathbb{H}^{n})$, one has
\begin{eqnarray*}
\|f\|_{L^{2}(\mathbb{H}^{n})}^{2}=\int\limits_{\mathbb{R}}
\|\widehat{f}(\lambda)\|_{\mathcal{B}_{2}}^{2}|\lambda|^{n}d\lambda,\label{2peq26}
\end{eqnarray*}
which is known as the Plancherel formula for the group Fourier transform on $\mathbb{H}^{n}$.
Define $f^\lambda(z)=\ds\int\limits_{\mathbb{R}}f(z,t)e^{2\pi i\lambda t}dt$ to be the inverse Fourier transform of $f$ in the $t$-variable. Thus $\widehat{f}(\lambda)=\ds\int\limits_{\mathbb{C}^n}f^\lambda(z)\pi_\lambda(z,0)dz$. Hence it becomes natural to consider the operator of the form
$$W_\lambda(g)=\int\limits_{\mathbb{C}^n}g(z)\pi_\lambda(z,0)dz,$$
for $g\in L^{1}(\mathbb{C}^{n})$. Then $\widehat{f}(\lambda)$ can be written as $\widehat{f}(\lambda)=W_{\lambda}(f^{\lambda})$. The transform $W_{\lambda}(g)$ of $g\in L^{1}(\mathbb{C}^{n})$ can be explicitly written as
\begin{eqnarray*}
W_{\lambda}(g)\varphi (\xi)=\int\limits_{\mathbb{C}^{n}}  g(z)e^{2\pi i\lambda (x.\xi+\frac{1}{2}x.y)}\varphi(\xi+y)dz,~~ \varphi\in L^{2}(\mathbb {R}^{n}),~~ z=x+iy,
\end{eqnarray*}
which maps $L^{1}(\mathbb {C}^{n})$ into $\mathcal{B}(L^{2}(\mathbb {R}^{n}))$. The transform $W_{\lambda}(g)$ is an integral operator with kernel $K_{g}^{\lambda}(\xi,\eta)$ given by
\begin{eqnarray*}
\int\limits_{\mathbb{R}^{n}}g(x,\eta-\xi)e^{i\pi\lambda x.(\xi+\eta)}dx.
\end{eqnarray*}

The convolution of $f,g\in L^{1}(\mathbb{H}^{n})$ is defined by
$$f\ast g(z,t)=\int\limits_{\mathbb{H}^{n}}f((z,t).(w,s)^{-1})g(w,s)dwds.$$
Further $L^{1}(\mathbb{H}^{n})$ turns out to be a non-commutative Banach algebra under convolution. The $\lambda$-twisted convolution of $F,G\in L^{1}(\mathbb{C}^{n})$ is defined to be
\begin{eqnarray*}
F\ast_{\lambda}G(z)=\int\limits_{\mathbb{C}^{n}}F(z-w)G(w)e^{\pi i\lambda Im(z.\overline{w})}dw.
\end{eqnarray*}
Now if $f,g\in L^{1}(\mathbb{H}^{n})$, then $f^{\lambda},g^{\lambda}\in L^{1}(\mathbb{C}^{n})$ and $(f\ast g)^{\lambda}=f^{\lambda}\ast_{\lambda} g^{\lambda}$. We refer to Thangavelu \cite{th} for further details on $\mathbb{H}^{n}$.

We organize the paper as follows. In section 2, we list properties of $\lambda$-twisted translations and introduce the map $\mathscr{T}$. In section 3, we obtain characterization of orthonormality of $E(\mathscr{A})$ in terms of the collection $\{\mathscr{T}\varphi_{j}(\lambda):j\in\mathbb{N}\}$ for a.e. $\lambda\in(0,1]$. In section 4, we prove that if $E(\mathscr{A})$ is a Bessel sequence then $\{\mathscr{T}\varphi_{j}(\lambda):j\in\mathbb{N}\}$ is a Bessel sequence for a.e. $\lambda\in(0,1]$ with the same bound. However, the converse becomes true with an additional condition ``condition C", given in Definition \ref{2pdef3}. In section 5, we prove that $E(\mathscr{A})$ is a Parseval frame if and only if  $\{\mathscr{T}\varphi_{j}(\lambda):j\in\mathbb{N}\}$ is a Parseval frame for a.e. $\lambda\in(0,1]$ under condition C. We also mention similar characterization for $E(\mathscr{A})$ to be a frame and Riesz basis. Also we obtain a decomposition theorem for a shift invariant space V in $L^{2}(\mathbb{H}^{n})$. In section 6, we focus on principal shift-invariant space $V(\varphi)$ and obtain certain additional results.
\section{Preliminaries}
Let $\mathcal{H}$ be a separable Hilbert space.
\begin{defn}
A sequence $\{f_{k}:{k \in\mathbb {Z}}\}$ in $\mathcal{H}$ is called a Bessel sequence for $\mathcal{H}$ if there exists a constant $B > 0$ such that
\begin{eqnarray*}
\sum \limits_{k \in\mathbb {Z}}|\langle f,f_{k} \rangle |^{2}\leq B\|f\|^{2},\hspace{4mm}\forall f\in \mathcal{H}.
\end{eqnarray*}
\end{defn}
\begin{defn}
A sequence $\{f_{k}:{k \in\mathbb {Z}}\}$ in $\mathcal{H}$ is called a frame for $\mathcal{H}$ if there exist two constants $A,B > 0$ such that
\begin{eqnarray*}
A\|f\|^{2}\leq\ds\sum \limits_{k \in\mathbb {Z}}|\langle f,f_{k} \rangle |^{2}\leq B\|f\|^{2},\hspace{4mm}\forall f\in \mathcal{H}.
\end{eqnarray*}
In particular, if $A=B=1$, then $\{f_{k}:{k \in\mathbb {Z}}\}$ is called a Parseval frame.

Let $S:\mathcal{H}\longrightarrow \mathcal{H}$ be defined by $Sf:= \sum\limits_{k\in\mathbb {Z}}\langle f,f_{k}\rangle f_{k}$, associated with the frame $\{f_{k}\}$. Then $S$ is called a frame operator, which is bounded, invertible, self-adjoint and positive. Further, $\{S^{-1}f_{k}:{k \in\mathbb {Z}}\}$ is also a frame for  $\mathcal{H}$ and is called the canonical dual frame of $\{f_{k}:{k \in\mathbb {Z}}\}$.
\end{defn}
\begin{defn}
Let V be a closed subspace of $\mathcal{H}$. Let $\{f_{k}:{k \in\mathbb {Z}}\}$ and $\{g_{k}:{k \in\mathbb {Z}}\}$ be two Bessel sequences in $\mathcal{H}$. If
$f=\ds\sum\limits_{k\in\mathbb{Z}}\langle f, g_{k}\rangle f_{k}$
holds for all $f\in V$, then  $\{g_{k}:{k \in\mathbb {Z}}\}$ is called an oblique dual of $\{f_{k}:{k \in\mathbb {Z}}\}$.
\end{defn}
\begin{defn}\label{2pdef1}
A Riesz basis for $\mathcal{H}$ is a family of the form $\{Ue_{k}:k\in\mathbb {Z}\}$, where $\{e_{k}:k\in\mathbb {Z}\}$ is an orthonormal basis for $\mathcal{H}$ and $U:\mathcal{H}\longrightarrow \mathcal{H}$ is a bounded invertible operator. Equivalently, a Riesz basis is a sequence  $\{f_{k}:k\in\mathbb {Z}\}$ which is complete in $\mathcal{H}$ and there exit constants $A, B> 0$ such that for every finite sequence of complex numbers $\{c_{k}\}$, one has
\begin{eqnarray*}\label{def2}
A\sum\limits_{k}\left|c_{k}\right|^{2}\leq \left\|\sum\limits_{k}c_{k}f_{k}\right\|^{2}\leq B\sum\limits_{k}\left|c_{k}\right|^{2}.
\end{eqnarray*}
\end{defn}
For a detailed study of frames, we refer to \cite{chr2} and \cite{heil}.
\begin{defn}
Let $f\in L^{2}(\mathbb{R}^{2n})$ and $\lambda\in\mathbb{R}^{\star}$. For $(k,l)\in \mathbb{Z}^{2n}$, $\lambda$-twisted translation of $f$, denoted by $(T_{(k,l)}^{t})^{\lambda}f$, is defined to be
\begin{eqnarray*}
(T_{(k,l)}^{t})^{\lambda}f(x,y)=e^{\pi i\lambda(x.l-y.k)}f (x-k,y-l),~~(x,y)\in\mathbb {R}^{2n}.
\end{eqnarray*}
\end{defn}
\subsection{Properties of $\lambda$-twisted translation}
\begin{enumerate}
\item The adjoint $((T_{(k,l)}^{t})^{\lambda})^{*}$ of $(T_{(k,l)}^{t})^{\lambda}$ is $(T_{(-k,-l)}^{t})^{\lambda}$.
\item $(T_{(k_{1}, l_{1})}^{t})^{\lambda}(T_{(k_{2}, l_{2})}^{t})^{\lambda}=
e^{-\pi i\lambda(k_{1}.l_{2}-k_{2}.l_{1})}(T_{(k_{1}+k_{2}, l_{1}+l_{2})}^{t})^{\lambda}$.
\item $(T_{(k,l)}^{t})^{\lambda}$ is a unitary operator on $L^{2}(\mathbb{R}^{2n})$
for all $(k,l)\in \mathbb{Z}^{2n}$.
\item The operator $W_{\lambda}$ satisfies $W_{\lambda}((T_{(k,l)}^{t})^{\lambda}f)=\pi_{\lambda}(k,l) W_{\lambda}(f)$.
\item Let $f\in L^{2}(\mathbb{C}^{n})$. Then the kernel of  $W_{\lambda}((T_{(k,l)}^{t})^{\lambda}f)$ satisfies the following relation.
    \begin{eqnarray*}
    K_{(T_{(k,l)}^{t})^{\lambda}f}^{\lambda}(\xi,\eta)=e^{\pi i\lambda(2\xi+l).k} K_{f}^{\lambda}(\xi+l,\eta).\label{2peq3}
    \end{eqnarray*}

\item Let $\varphi\in L^{2}(\mathbb{H}^{n})$. Then
\begin{eqnarray*}
 (L_{(2k,l,m)}\varphi)^{\lambda}= e^{2\pi im\lambda}(T_{(2k,l)}^{t})^{\lambda}\varphi^{\lambda}.\label{2peq4}
 \end{eqnarray*}
\end{enumerate}

Let $\mathscr{T}:L^{2}(\mathbb{H}^{n})\rightarrow L^{2}((0,1],\ell^{2}(\mathbb{Z},\mathcal{B}_{2}))$ be a map defined for $f\in L^{2}(\mathbb{H}^{n})$ by
$$\mathscr{T}f: (0,1]\rightarrow\ell^{2}(\mathbb{Z},\mathcal{B}_{2}),~~\mathscr{T}f(\lambda)=
\{\widehat{f}(\lambda+r)|\lambda+r|^{\frac{n}{2}}\}_{r\in\mathbb{Z}}.$$
Then $\mathscr{T}$ is an isometric isomorphism between $L^{2}(\mathbb{H}^{n})$ and $L^{2}((0,1],\ell^{2}(\mathbb{Z},\mathcal{B}_{2}))$.
\section{Orthonormality of $E(\mathscr{A})$ in $L^{2}(\mathbb{H}^{n})$}
In this section, we shall study orthonormality of $E(\mathscr{A})$ in $L^{2}(\mathbb{H}^{n})$. In fact we show that orthonormality of  $E(\mathscr{A})$ implies  orthonormality of $\{\mathscr{T}\varphi_{j}(\lambda):j\in\mathbb{N}\}$ in $\ell^{2}(\mathbb{Z},\mathcal{B}_{2})$ for a.e. $\lambda\in(0,1]$. But for the converse, orthonormality of $\{\mathscr{T}\varphi_{j}(\lambda):j\in\mathbb{N}\}$  for a.e. $\lambda\in(0,1]$ together with certain additional condition (Definition \ref{2pdef3}) implies the orthonormality of $E(\mathscr{A})$. We also show that the additional condition can not be dropped. (See Example \ref{exam1}).

Using the kernel of $W_{\lambda}(\varphi^{\lambda})$, for $\varphi\in L^{2}(\mathbb{H}^{n})$, we introduce the following function.
For each $k,l\in\mathbb{Z}^{n}$, we define
\begin{eqnarray*}
&&G_{k,l}^{\varphi}(\lambda)=\sum\limits_{r\in\mathbb{Z}}\sum\limits_{s\in\mathbb{Z}^{n}}
\int\limits_{\mathbb{T}^{n}}\int\limits_{\mathbb{R}^{n}}
K^{\lambda+r}_{\varphi^{\lambda+r}}(\xi+s,\eta)
\overline{K^{\lambda+r}_{\varphi^{\lambda+r}}(\xi+s+l,\eta)}\\
&&\hspace{3cm}\times e^{-2\pi i(\lambda+r)k.(l+2\xi)} e^{-4\pi i\lambda k.s} d\eta d\xi|\lambda+r|^{n},~~~\lambda\in (0,1].
\end{eqnarray*}

\begin{defn}\label{2pdef3}
A function $\varphi\in L^{2}(\mathbb{H}^{n})$ is said to satisfy condition C if $G_{k,l}^{\varphi}(\lambda)=0$ a.e. $\lambda\in (0,1]$,~for all $(k,l)\in\mathbb{Z}^{2n}\backslash\{(0,0)\}$.
\end{defn}
Now, we shall prove the following result.
\begin{theorem}\label{2pth1}
If $\varphi\in L^{2}(\mathbb{H}^{n})$, then $\{L_{(2k,l,m)}\varphi: (k,l,m)\in\mathbb{Z}^{2n+1}\}$ is an orthonormal system in $L^{2}(\mathbb{H}^{n})$ if and only if $G_{0,0}^{\varphi}(\lambda)=1$ a.e. $\lambda\in (0,1]$ and $\varphi$ satisfies condition C.
\end{theorem}
In order to prove this theorem, we prove the following Lemmas.
\begin{lem}\label{2plem1}
Let $\varphi\in L^{2}(\mathbb{H}^{n})$. If $\{L_{(2k,l,m)}\varphi: (k,l,m)\in \mathbb{Z}^{2n+1}\}$ forms an orthonormal system in ${L^{2}(\mathbb{H}^{n})}$, then $G_{0,0}^{\varphi}(\lambda)=1$ a.e. $\lambda\in (0,1]$.
\end{lem}
\begin{lem}\label{2plem2}
Let $\varphi\in L^{2}(\mathbb{H}^{n})$. Assume that $G_{0,0}^{\varphi}(\lambda)=1$ a.e. $\lambda\in (0,1]$. If, in addition $\varphi$ satisfies condition C, then $\{L_{(2k,l,m)}\varphi: (k,l,m)\in \mathbb{Z}^{2n+1}\}$ forms an orthonormal system in ${L^{2}(\mathbb{H}^{n})}$.
\end{lem}
\begin{lem}\label{2plem3}
Let $\varphi\in L^{2}(\mathbb{H}^{n})$. If $\{L_{(2k,l,m)}\varphi: (k,l,m)\in \mathbb{Z}^{2n+1}\}$ is an orthonormal system in ${L^{2}(\mathbb{H}^{n})}$, then $\varphi$ satisfies condition C.
\end{lem}
Combining Lemma~\ref{2plem1}, Lemma~\ref{2plem2} and Lemma~\ref{2plem3}, we get Theorem \ref{2pth1}.

Now, we shall proceed towards proving our main result.
\begin{lem}\label{2plem4}
For $\varphi_{i},\varphi_{j}\in\mathscr{A}$ with $i\neq j$ and $k_{1},k_{2},l_{1},l_{2}\in\mathbb{Z}^{n}$, one has the following identity.
\begin{eqnarray*}
\sum\limits_{r\in\mathbb{Z}}\langle \widehat{L_{(2k_{1},l_{1},0)}\varphi_{i}}(\lambda+r),\widehat{L_{(2k_{2},l_{2},0)}\varphi_{j}}
(\lambda+r)\rangle_{\mathcal{B}_{2}}|\lambda+r|^{n}=0~a.e.~\lambda\in (0,1].
\end{eqnarray*}
In particular, when $k_{1}=k_{2}=l_{1}=l_{2}=0$, one has
\begin{eqnarray}
\sum\limits_{r\in\mathbb{Z}}\langle\widehat{\varphi_{i}}(\lambda+r),
\widehat{\varphi_{j}}(\lambda+r)\rangle_{\mathcal{B}_{2}}|\lambda+r|^{n}=0 ~a.e.~\lambda\in (0,1].\label{2peq6}
\end{eqnarray}
\end{lem}
Using Theorem \ref{2pth1} and Lemma \ref{2plem4}, we get
\begin{theorem}\label{2pth2}
The collection $E(\mathscr{A})$ is an orthonormal system in $L^{2}(\mathbb{H}^{n})$ if and only if $\{\mathscr{T}\varphi_{j}(\lambda):j\in\mathbb{N}\}$ is an orthonormal system in $\ell^{2}(\mathbb{Z},\mathcal{B}_{2})$ for a.e. $\lambda\in (0,1]$ and each $\varphi_{j}$ satisfies condition C.
\end{theorem}
\begin{rem}\label{2prem1}
When $\mathscr{A}=\{\varphi\}$, then Theorem \ref{2pth2} boils down to Theorem \ref{2pth1}.
\end{rem}
\begin{rem}\label{2prem2}
The equivalent statement of Theorem \ref{2pth1} can be written in terms of bracket map, which is defined to be
\begin{eqnarray*}
[\varphi,\psi](\lambda)=\sum\limits_{r\in\mathbb{Z}}\langle\widehat{\varphi}
(\lambda+r),\widehat{\psi}(\lambda+r)\rangle_{\mathcal{B}_{2}}|\lambda+r|^{n},
~\lambda\in (0,1],
\end{eqnarray*}
for $\varphi,\psi\in L^{2}(\mathbb{H}^{n})$. Using this bracket map, the statement of Theorem \ref{2pth1} can be rephrased as follows:

If $\varphi\in L^{2}(\mathbb{H}^{n})$, then $\{L_{(2k,l,m)}\varphi: (k,l,m)\in\mathbb{Z}^{2n+1}\}$ is an orthonormal system in $L^{2}(\mathbb{H}^{n})$ if and only if $[\varphi, L_{(2k,l,0)}\varphi](\lambda)=\delta_{k,l,0}$ a.e. $\lambda\in (0,1],\forall (k,l)\in\mathbb{Z}^{2n}$. This follows from the fact that $[\varphi, L_{(2k,l,0)}\varphi](\lambda)=G_{k,l}^{\varphi}(\lambda)$ for $\lambda\in (0,1]$. 
In this connection we refer to \cite{bhm}, wherein a similar result is proved for the polarized Heisenberg group $\mathbb{H}^{n}$. \textbf{We wish to emphasize that the proof given in this paper is totally different from the proof given in \cite{bhm}.}
\end{rem}
\begin{example}\label{exam1}
In Theorem \ref{2pth1}, condition C of $\varphi$ can not be dropped. We shall give an example of $\varphi$ which satisfies $G_{0,0}^{\varphi}(\lambda)=1$ a.e. $\lambda\in (0,1]$ but does not satisfy condition C. In this case we shall see that $\{L_{(2k,l,m)}\varphi: (k,l,m)\in\mathbb{Z}^{2n+1}\}$ is not an orthonormal system.
\end{example}
Let $\varphi=\frac{1}{3}\chi_{[0,3]\times [0,3]\times [0,1]}$. Then $\{L_{(2k,l,m)}\varphi: (k,l,m)\in\mathbb{Z}^{3}\}$ is not an orthonormal system in $L^{2}(\mathbb{H}^{1})$ but $G_{0,0}^{\varphi}(\lambda)=1$ a.e. $\lambda\in (0,1]$. 
\begin{example}
Condition C is weaker than the orthonormality of $\{L_{(2k,l,m)}\varphi: (k,l,m)$ $\in\mathbb{Z}^{2n+1}\}$. We shall illustrate this with an example.
\end{example}
Let $\varphi=\chi_{[0,2]\times [0,1]\times [0,2]}$. Then $\langle\varphi, L_{(0,0,1)}{\varphi}\rangle=2\neq 0$, from which it follows that $\{L_{(2k,l,m)}\varphi:{(k,l,m)\in\mathbb{Z}^{3}}$\} does not form an orthonormal system in $ L^{2}(\mathbb{H}^{1})$ but $\varphi$ satisfies condition C.
\section{$E(\mathscr{A})$ as a Bessel sequence in $L^{2}(\mathbb{H}^{n})$}
\begin{theorem}\label{2pth6}
Suppose $E(\mathscr{A})$ is a Bessel sequence in $L^{2}(\mathbb{H}^{n})$ with bound B. Then $\{\mathscr{T}\varphi_{j}(\lambda):j\in\mathbb{N}\}$ is a Bessel sequence in $\ell^{2}(\mathbb{Z},\mathcal{B}_{2})$ with bound B for a.e. $\lambda\in (0,1]$.

Conversely, suppose $\{\mathscr{T}\varphi_{j}(\lambda):j\in\mathbb{N}\}$ is a Bessel sequence in $\ell^{2}(\mathbb{Z},\mathcal{B}_{2})$ with bound B for a.e. $\lambda\in (0,1]$. If, in addition each $\varphi_{j}$ satisfies condition C, then $E(\mathscr{A})$ is a Bessel sequence in $L^{2}(\mathbb{H}^{n})$ with bound B.
\end{theorem}

\begin{pf}
Let $E(\mathscr{A})$ be a Bessel sequence in $L^{2}(\mathbb{H}^{n})$ with bound B. Then
\begin{eqnarray}
\sum\limits_{j\in\mathbb{N}}\sum\limits_{(k,l,m)\in\mathbb{Z}^{2n+1}} |\langle f,L_{(2k,l,m)}\varphi_{j}\rangle_{L^{2}(\mathbb{H}^{n})}|^{2}\leq B\|f\|_{L^{2}(\mathbb{H}^{n})}^{2},~\forall~f\in L^{2}(\mathbb{H}^{n}).\label{2peq9}
\end{eqnarray}
But
\begin{eqnarray*}
&&\sum\limits_{j\in\mathbb{N}}\sum\limits_{(k,l,m)\in\mathbb{Z}^{2n+1}} |\langle f,L_{(2k,l,m)}\varphi_{j}\rangle|^{2}\\
&\geq&\sum\limits_{j\in\mathbb{N}}\sum\limits_{m\in\mathbb{Z}}|\langle f,L_{(0,0,m)}\varphi_{j}\rangle|^{2}\\
&=&\sum\limits_{j\in\mathbb{N}}\sum\limits_{m\in\mathbb{Z}} \bigg|\int\limits_{\mathbb{R}}\langle\widehat{f}(\lambda), \widehat{L_{(0,0,m)}\varphi_{j}}(\lambda)\rangle_{\mathcal{B}_{2}}
|\lambda|^{n}d\lambda\bigg|^{2}\\
&=&\sum\limits_{j\in\mathbb{N}}\sum\limits_{m\in\mathbb{Z}} \bigg|\int\limits_{\mathbb{R}}\langle\widehat{f}(\lambda),
\widehat{\varphi_{j}}(\lambda)\rangle_{\mathcal{B}_{2}}|\lambda|^{n} e^{-2\pi im\lambda}d\lambda\bigg|^{2},
\end{eqnarray*}
using the fact that $\widehat{L_{(x,y,t)}\varphi}(\lambda)=e^{2\pi i\lambda t}\widehat{L_{(x,y,0)}\varphi}(\lambda),\forall~(x,y,t)\in\mathbb{H}^{n}$. Hence
\begin{eqnarray*}
&&\sum\limits_{j\in\mathbb{N}}\sum\limits_{(k,l,m)\in\mathbb{Z}^{2n+1}} |\langle f,L_{(2k,l,m)}\varphi_{j}\rangle|^{2}\\
&\geq&\sum\limits_{j\in\mathbb{N}}\sum\limits_{m\in\mathbb{Z}} \bigg|\int\limits_{0}^{1}\sum\limits_{r\in\mathbb{Z}}\langle\widehat{f}(\lambda+r),
\widehat{\varphi_{j}}(\lambda+r)
\rangle_{\mathcal{B}_{2}}|\lambda+r|^{n} e^{-2\pi im\lambda}d\lambda\bigg|^{2}\\
&=&\sum\limits_{j\in\mathbb{N}}\int\limits_{0}^{1}\bigg|\sum\limits_{r\in\mathbb{Z}}
\langle\widehat{f}(\lambda+r),
\widehat{\varphi_{j}}(\lambda+r)
\rangle_{\mathcal{B}_{2}}|\lambda+r|^{n}\bigg|^{2}d\lambda,
\end{eqnarray*}
using Parseval's formula for the Fourier coefficient with respect to the $m$ variable.
Thus
\begin{eqnarray*}
\sum\limits_{j\in\mathbb{N}}\sum\limits_{(k,l,m)\in\mathbb{Z}^{2n+1}} |\langle f,L_{(2k,l,m)}\varphi_{j}\rangle|^{2}
\geq\sum\limits_{j\in\mathbb{N}}\int\limits_{0}^{1}\bigg|\sum\limits_{r\in\mathbb{Z}}
\langle\widehat{f}(\lambda+r),
\widehat{\varphi_{j}}(\lambda+r)
\rangle_{\mathcal{B}_{2}}|\lambda+r|^{n}\bigg|^{2}d\lambda.
\end{eqnarray*}
Hence it follows from (\ref{2peq9}) that
\begin{eqnarray*}
\sum\limits_{j\in\mathbb{N}}\int\limits_{0}^{1}\bigg|\sum\limits_{r\in\mathbb{Z}}
\langle\widehat{f}(\lambda+r),
\widehat{\varphi_{j}}(\lambda+r)
\rangle_{\mathcal{B}_{2}}|\lambda+r|^{n}\bigg|^{2}d\lambda\leq B\|f\|^{2}.
\end{eqnarray*}
But
\begin{eqnarray*}
\|f\|_{L^{2}(\mathbb{H}^{n})}^{2}=\|\mathscr{T}f\|_{L^{2}((0,1],\ell^{2}(\mathbb{Z},
\mathcal{B}_{2}))}^{2}=\int\limits_{0}^{1}\sum\limits_{r\in\mathbb{Z}}
\|\widehat{f}(\lambda+r)\|_{\mathcal{B}_{2}}^{2}|\lambda+r|^{n}d\lambda.
\end{eqnarray*}
Hence
\begin{eqnarray}
&&\hspace{-4cm}\sum\limits_{j\in\mathbb{N}}\int\limits_{0}^{1}\bigg
|\sum\limits_{r\in\mathbb{Z}}\langle\widehat{f}(\lambda+r),
\widehat{\varphi_{j}}(\lambda+r)
\rangle_{\mathcal{B}_{2}}|\lambda+r|^{n}\bigg|^{2}d\lambda\nonumber\\
&\leq& B\int\limits_{0}^{1}\sum\limits_{r\in\mathbb{Z}}
\|\hat{f}(\lambda+r)\|_{\mathcal{B}_{2}}^{2}|\lambda+r|^{n}d\lambda,~\forall f\in L^{2}(\mathbb{H}^{n}).\label{2peq10}
\end{eqnarray}
Now, our claim is to show that for a.e. $\lambda\in (0,1]$,
\begin{eqnarray}
G_{0,0}^{\varphi_{j}}(\lambda)\leq B,~\forall j\in\mathbb{N}.\label{2peq11}
\end{eqnarray}
If it were not true, then there would exists some $j_{0}\in\mathbb{N}$ such that $M=\{\lambda\in (0,1]:\sum\limits_{r\in\mathbb{Z}}\|\widehat{\varphi_{j_{0}}}(\lambda+r)\|
_{\mathcal{B}_{2}}^{2}|\lambda+r|^{n}> B\}$ has strictly positive measure. Let $g(\lambda)=\chi_{M}(\lambda)$ on (0,1]and extend $g$ on $\mathbb{R}$ by defining $g(\lambda)=g(\lambda+1)$ for all $\lambda\in\mathbb{R}$. Then there exists  $f_{j_{0}}\in L^{2}(\mathbb{H}^{n})$ such that $\widehat{f_{j_{0}}}(\lambda)=\chi_{M}(\lambda)\widehat{\varphi_{j_{0}}}(\lambda)$. Now
\begin{eqnarray*}
&&\int\limits_{0}^{1}\sum\limits_{j\in\mathbb{N}}\bigg|\sum\limits_{r\in\mathbb{Z}}
\langle\widehat{f_{j_{0}}}(\lambda+r),\widehat{\varphi_{j}}(\lambda+r)\rangle
_{\mathcal{B}_{2}}|\lambda+r|^{n}\bigg|^{2}d\lambda\\
&\geq&\int\limits_{0}^{1}\bigg|\sum\limits_{r\in\mathbb{Z}}
\langle\widehat{f_{j_{0}}}(\lambda+r),\widehat{\varphi_{j_{0}}}(\lambda+r)\rangle
_{\mathcal{B}_{2}}|\lambda+r|^{n}\bigg|^{2}d\lambda\\
&=&\int\limits_{0}^{1}\bigg|\sum\limits_{r\in\mathbb{Z}}
\langle\chi_{M}(\lambda)\widehat{\varphi_{j_{0}}}(\lambda+r),\widehat{\varphi_{j_{0}}}
(\lambda+r)\rangle_{\mathcal{B}_{2}}
|\lambda+r|^{n}\bigg|^{2}d\lambda\\
&=&\int\limits_{M}\bigg(\sum\limits_{r\in\mathbb{Z}}\|\widehat{\varphi_{j_{0}}}(\lambda+r)\|
_{\mathcal{B}_{2}}^{2}|\lambda+r|^{n}\bigg)^{2}d\lambda\\
&>& B\int\limits_{M}\sum\limits_{r\in\mathbb{Z}}\|\widehat{\varphi_{j_{0}}}(\lambda+r)\|
_{\mathcal{B}_{2}}^{2}|\lambda+r|^{n}d\lambda\\
&=& B\int\limits_{0}^{1}\sum\limits_{r\in\mathbb{Z}}\|\chi_{M}(\lambda)
\widehat{\varphi_{j_{0}}}(\lambda+r)\|
_{\mathcal{B}_{2}}^{2}|\lambda+r|^{n}d\lambda\\
&=&B\int\limits_{0}^{1}\sum\limits_{r\in\mathbb{Z}}\|\widehat{f_{j_{0}}}(\lambda+r)\|
_{\mathcal{B}_{2}}^{2}|\lambda+r|^{n}d\lambda,
\end{eqnarray*} which is a contradiction to (\ref{2peq10}). This shows that M has measure zero, from which our claim (\ref{2peq11}) follows. Now let $\{c_{j}:j\in\mathbb{N}\}$ be a finite sequence. It is enough to show that for a.e. $\lambda\in (0,1]$, $\bigg\|\sum\limits_{j\in\mathcal{F}}c_{j}\mathscr{T}\varphi_{j}(\lambda)\bigg\|
_{\ell^{2}(\mathbb{Z},\mathcal{B}_{2})}^{2}\leq B\sum\limits_{j\in\mathcal{F}}|c_{j}|^{2}$, where $\mathcal{F}$ denotes a finite set. Now, from the definition of $G_{k,l}^{\varphi_{j}}$, we have
\begin{eqnarray}
G_{0,0}^{\varphi_{j}}(\lambda)&=&\sum\limits_{r\in\mathbb{Z}}\sum\limits_
{s\in\mathbb{Z}^{n}}\int\limits_{\mathbb{T}^{n}}\int\limits_{\mathbb{R}^{n}}
|K^{\lambda+r}_{\varphi_{j}^{\lambda+r}}(\xi+s,\eta)|^{2}d\eta d\xi|\lambda+r|^{n}\nonumber\\
&=&\sum\limits_{r\in\mathbb{Z}}\int\limits_{\mathbb{R}^{n}}\int\limits_
{\mathbb{R}^{n}}|K^{\lambda+r}_{\varphi_{j}^{\lambda+r}}(\xi,\eta)|^{2}d\eta
d\xi|\lambda+r|^{n}\nonumber\\
&=&\sum\limits_{r\in\mathbb{Z}}\|W_{\lambda+r}(\varphi_{j}^{\lambda+r})\|
_{\mathcal{B}_{2}}^{2}|\lambda+r|^{n}\nonumber\\
&=&\sum\limits_{r\in\mathbb{Z}}\|\widehat{\varphi_{j}}(\lambda+r)\|_
{\mathcal{B}_{2}}^{2}|\lambda+r|^{n}\nonumber\\
&=&\|\mathscr{T}\varphi_{j}(\lambda)\|_{\ell^{2}(\mathbb{Z},\mathcal{B}_{2})}^{2}.
\label{2peq34}
\end{eqnarray}
Consider for $\lambda\in (0,1]$,
\begin{eqnarray*}
\bigg\|\sum\limits_{j\in\mathcal{F}}c_{j}\mathscr{T}\varphi_{j}(\lambda)\bigg\|
_{\ell^{2}(\mathbb{Z},\mathcal{B}_{2})}^{2}
=\sum\limits_{j}|c_{j}|^{2}\|\mathscr{T}
\varphi_{j}(\lambda)\|^{2}+\sum\limits_{j_{1}\neq j_{2}}c_{j_{1}}\overline{c_{j_{2}}}\langle\mathscr{T}\varphi_{j_{1}}(\lambda),
\mathscr{T}\varphi_{j_{2}}(\lambda)\rangle .
\end{eqnarray*}
But for $j_{1}\neq j_{2}$,
\begin{eqnarray*}
\langle\mathscr{T}\varphi_{j_{1}}(\lambda),
\mathscr{T}\varphi_{j_{2}}(\lambda)\rangle
=\sum\limits_{r\in\mathbb{Z}}\langle\widehat{\varphi_{j_{1}}}(\lambda+r),
\widehat{\varphi_{j_{2}}}(\lambda+r)\rangle|\lambda+r|^{n}=0~a.e.~\lambda\in (0,1],
\end{eqnarray*}
using (\ref{2peq6}).
Thus for a.e. $\lambda\in (0,1]$,
\begin{eqnarray}
\bigg\|\sum\limits_{j\in\mathcal{F}}c_{j}\mathscr{T}\varphi_{j}(\lambda)\bigg\|
_{\ell^{2}(\mathbb{Z},\mathcal{B}_{2})}^{2}
=\sum\limits_{j}|c_{j}|^{2}\|\mathscr{T}
\varphi_{j}(\lambda)\|^{2}
&=&\sum\limits_{j}|c_{j}|^{2}G_{0,0}^{\varphi_{j}}(\lambda)\label{2peq8}\\
&\leq& B\sum\limits_{j\in\mathcal{F}} |c_{j}|^{2},\nonumber
\end{eqnarray}
using (\ref{2peq34}) and (\ref{2peq11}).
Hence $\{\mathscr{T}\varphi_{j}(\lambda):j\in\mathbb{N}\}$ is a Bessel sequence in $\ell^{2}(\mathbb{Z},\mathcal{B}_{2})$ with bound B for a.e. $\lambda\in (0,1]$.

Conversely, suppose $\{\mathscr{T}\varphi_{j}(\lambda):j\in\mathbb{N}\}$ is a Bessel sequence in $\ell^{2}(\mathbb{Z},\mathcal{B}_{2})$ with bound B for a.e. $\lambda\in (0,1]$. This means
\begin{eqnarray}
\sum\limits_{j\in\mathbb{N}}|\langle a,\mathscr{T}\varphi_{j}(\lambda)\rangle|^{2}\leq B\|a\|^{2}, \forall~ a\in\ell^{2}(\mathbb{Z},\mathcal{B}_{2}).\label{2peq12}
\end{eqnarray}
Now put $a=\mathscr{T}f(\lambda)$ in (\ref{2peq12}), for $f\in L^{2}(\mathbb{H}^{n})$. Then integrating (\ref{2peq12}) over [0,1], we get (\ref{2peq10}). This in turn will lead to (\ref{2peq11}), as discussed earlier.
In order to prove our result, it is enough to show that
$\bigg\|\sum\limits_{(k,l,m,j)\in\mathcal{F}}c_{k,l,m,j}L_{(2 k,l,m)}\varphi_{j}\bigg\|_{L^{2}(\mathbb{H}^{n})}^{2}\leq B\sum\limits_{(k,l,m,j)\in\mathcal{F}}|c_{k,l,m,j}|^{2}$, for a finite set $\mathcal{F}$. Now, consider
\begin{eqnarray*}
\|\widehat{L_{(2k,l,0)}\varphi}(\lambda)\|_{\mathcal{B}_{2}} =\|\pi_{\lambda}(2k,l,0)\widehat{\varphi}(\lambda)\|_{\mathcal{B}_{2}}
&=& \|K_{\pi_{\lambda}(2k,l,0)\widehat{\varphi}(\lambda)}\|_{L^{2}(\mathbb{C}^{n})}\\
&=&\|K_{\pi_{\lambda}(2k,l,0)W_{\lambda}(\varphi^{\lambda})}\|_{L^{2}
(\mathbb{C}^{n})}.
\end{eqnarray*}
It can be easily shown that
\begin{eqnarray*}
K_{\pi_{\lambda}(2k,l,0)W_{\lambda}(\varphi^{\lambda})}(\xi,\eta)=e^{2\pi i\lambda(2k.\xi+k.l)}K_{\varphi^{\lambda}}^{\lambda}(\xi+l,\eta).
\end{eqnarray*}
Thus
\begin{eqnarray}
\|\widehat{L_{(2k,l,0)}\varphi}(\lambda)\|_{\mathcal{B}_{2}}^{2}
=\int\limits_{\mathbb{R}^{2n}}|K^{\lambda}_{\varphi^{\lambda}}(\xi+l,\eta)|^{2}
d\xi d\eta
=\|K^{\lambda}_{\varphi^{\lambda}}\|_{L^{2}(\mathbb{C}^{n})}^{2}
=\|W^{\lambda}(\varphi^{\lambda})\|_{\mathcal{B}_{2}}^{2}
=\|\widehat{\varphi}(\lambda)\|_{\mathcal{B}_{2}}^{2}.\nonumber\\ \label{2peq14}
\end{eqnarray}
Now, consider
\begin{eqnarray}
&&\bigg\|\sum\limits_{(k,l,m,j)\in\mathcal{F}}c_{k,l,m,j}L_{(2    k,l,m)}\varphi_{j}\bigg\|_{L^{2}(\mathbb{H}^{n})}^{2}\nonumber\\
&=&\int\limits_{\mathbb{R}}\bigg\|\sum\limits_{k,l,m,j}c_{k,l,m,j}
\widehat{L_{(2k,l,m)}\varphi_{j}(\lambda)}\bigg\|_{\mathcal{B}_{2}}^{2}
|\lambda|^{n}d\lambda\nonumber\\
&=&\int\limits_{0}^{1}\sum\limits_{r\in\mathbb{Z}}\bigg\|\sum\limits_{(k,l,m,j)}
c_{k,l,m,j}\widehat{L_{(2k,l,m)}\varphi_{j}}(\lambda+r)\bigg\|
_{\mathcal{B}_{2}}^{2}|\lambda+r|^{n}d\lambda\nonumber\\
&=&\int\limits_{0}^{1}\sum\limits_{r\in\mathbb{Z}}\bigg\|\sum\limits_{k,l,j}
\bigg(\sum\limits_{m}c_{k,l,m,j}e^{2\pi im\lambda}\bigg)\widehat{L_{(2k,l,0)}\varphi_{j}}(\lambda+r)\bigg\|
_{\mathcal{B}_{2}}^{2}|\lambda+r|^{n}d\lambda\nonumber\\
&=&\int\limits_{0}^{1}\sum\limits_{r\in\mathbb{Z}}\sum\limits_{k,l,j}\bigg\|
\bigg(\sum\limits_{m}c_{k,l,m,j}e^{2\pi im\lambda}\bigg)\widehat{L_{(2k,l,0)}\varphi_{j}}(\lambda+r)\bigg\|
_{\mathcal{B}_{2}}^{2}|\lambda+r|^{n}d\lambda\nonumber\\
&+&\int\limits_{0}^{1}\sum\limits_{r\in\mathbb{Z}}
\sum\limits_{(k_{1},l_{1},j_{1})\neq (k_{2},l_{2},j_{2})}\bigg\langle\bigg(\sum\limits_{m}c_{k_{1},l_{1},m,j_{1}}e^{2\pi i m\lambda}\bigg)\widehat{L_{(2k_{1},l_{1},0)}\varphi_{j_{1}}}(\lambda+r),\nonumber\\
&&\bigg(\sum\limits_{m}c_{k_{2},l_{2},m,j_{2}}e^{2\pi i m\lambda}\bigg)\widehat{L_{(2k_{2},l_{2},0)}\varphi_{j_{2}}}(\lambda+r)\bigg\rangle
_{\mathcal{B}_{2}}|\lambda+r|^{n}d\lambda.\label{2peq13}
\end{eqnarray}
Using (\ref{2peq14}) and (\ref{2peq11}), the first term on the right hand side of (\ref{2peq13}) becomes
\begin{eqnarray}
&&\int\limits_{0}^{1}\sum\limits_{r\in\mathbb{Z}}\sum\limits_{k,l,j}
\bigg|\sum\limits_{m}c_{k,l,m,j}e^{2\pi im\lambda}\bigg|^{2}\|\widehat{L_{(2k,l,0)}\varphi_{j}}(\lambda+r)\|
_{\mathcal{B}_{2}}^{2}|\lambda+r|^{n}d\lambda\nonumber\\
&=&\int\limits_{0}^{1}\sum\limits_{k,l,j}
\bigg|\sum\limits_{m}c_{k,l,m,j}e^{2\pi im\lambda}\bigg|^{2}\sum\limits_{r\in\mathbb{Z}}\|\widehat{\varphi_{j}}
(\lambda+r)\|_{\mathcal{B}_{2}}^{2}|\lambda+r|^{n}d\lambda\nonumber\\
&=&\int\limits_{0}^{1}\sum\limits_{k,l,j}
\bigg|\sum\limits_{m}c_{k,l,m,j}e^{2\pi im\lambda}\bigg|^{2}G_{0,0}^{\varphi_{j}}(\lambda)d\lambda\label{2peq15}\\
&\leq& B\sum\limits_{k,l,j}\int\limits_{0}^{1}\bigg|\sum\limits_{m}c_{k,l,m,j}e^{2\pi im\lambda}\bigg|^{2}d\lambda\nonumber\\
&=& B\sum\limits_{k,l,m,j}|c_{k,l,m,j}|^{2}\nonumber.
\end{eqnarray}
The second term on the right hand side of (\ref{2peq13}) is
\begin{eqnarray*}
&&\int\limits_{0}^{1}
\sum\limits_{(k_{1},l_{1},j_{1})\neq (k_{2},l_{2},j_{2})}\sum\limits_{m_{1},m_{2}}
c_{k_{1},l_{1},m_{1},j_{1}}\overline{c_{k_{2},l_{2},m_{2},j_{2}}}e^{2\pi i(m_{1}-m_{2})\lambda}\\
&&\times\sum\limits_{r\in\mathbb{Z}}\langle\widehat{L_{(2k_{1},l_{1},0)}\varphi_{j_{1}}}(\lambda+r),
\widehat{L_{(2k_{2},l_{2},0)}\varphi_{j_{2}}}(\lambda+r)\rangle
_{\mathcal{B}_{2}}|\lambda+r|^{n}d\lambda\\
&=&\int\limits_{0}^{1}\sum\limits_{(k_{1},l_{1})\neq (k_{2},l_{2})}\sum\limits_{j_{1}=j_{2}=j}\sum\limits_{m_{1},m_{2}}
c_{k_{1},l_{1},m_{1},j}\overline{c_{k_{2},l_{2},m_{2},j}}e^{2\pi i(m_{1}-m_{2})\lambda}\\
&&\times\sum\limits_{r\in\mathbb{Z}}\langle \widehat{L_{(2k_{1},l_{1},0)}\varphi_{j}}(\lambda+r),\widehat{L_{(2k_{2},l_{2},0)}
\varphi_{j}}(\lambda+r)\rangle_{\mathcal{B}_{2}}|\lambda+r|^{n}d\lambda\\
&+&\int\limits_{0}^{1}\sum\limits_{k_{1},k_{2},l_{1},l_{2}}\sum\limits_{j_{1}\neq j_{2}}\sum\limits_{m_{1},m_{2}}c_{k_{1},l_{1},m_{1},j_{1}}
\overline{c_{k_{2},l_{2},m_{2},j_{2}}}e^{2\pi i(m_{1}-m_{2})\lambda}\\
&&\times\sum\limits_{r\in\mathbb{Z}}\langle \widehat{L_{(2k_{1},l_{1},0)}\varphi_{j_{1}}}(\lambda+r),\widehat{L_{(2k_{2},l_{2},0)}
\varphi_{j_{2}}}(\lambda+r)\rangle_{\mathcal{B}_{2}}|\lambda+r|^{n}d\lambda.
\end{eqnarray*}
Here, since each $\varphi_{j}$ satisfies condition C, the first term is zero. Using Lemma \ref{2plem4}, the second term becomes zero. Thus
\begin{eqnarray*}
\bigg\|\sum\limits_{(k,l,m,j)\in\mathcal{F}}c_{k,l,m,j}L_{(2    k,l,m)}\varphi_{j}\bigg\|_{L^{2}(\mathbb{H}^{n})}^{2}\leq B
\sum\limits_{(k,l,m,j)\in\mathcal{F}}|c_{k,l,m,j}|^{2},
\end{eqnarray*}
proving that $E(\mathscr{A})$ is a Bessel sequence in $L^{2}(\mathbb{H}^{n})$ with bound B.
\end{pf}

\section{$E(\mathscr{A})$ as a Parseval frame}
Assume that each $\varphi_{j}$ satisfies condition C. Consider $f\in U(\mathscr{A})$ i.e.,\\ $f=\sum\limits_{(k^{\prime},l^{\prime},m^{\prime},j^{\prime})
\in\mathcal{F}}c_{k^{\prime},l^{\prime},m^{\prime},j^{\prime}} L_{(2k^{\prime},l^{\prime},m^{\prime})}\varphi_{j^{\prime}}$, where $\mathcal{F}$ is a finite set. Define $\fontdimen16\textfont2=5pt\rho_{j^{\prime}}(\lambda)=\{\fontdimen16\textfont2=5pt\rho_{k^{\prime},l^{\prime},j^{\prime}}
(\lambda)\}_{k^{\prime},l^{\prime}}$ for $\lambda\in (0,1]$, where
\begin{eqnarray*}
\fontdimen16\textfont2=5pt\rho_{k^{\prime},l^{\prime},j^{\prime}}
(\lambda)=\sum\limits_{m^{\prime}}c_{k^{\prime},l^{\prime},m^{\prime},j^{\prime}}
e^{2\pi im^{\prime}\lambda}\label{2peq16}
\end{eqnarray*}
and
\begin{eqnarray*}
\rho(\lambda)=\bigoplus\limits_{j^{\prime}}\fontdimen16\textfont2=5pt\rho_{j^{\prime}}(\lambda).\label{2peq17}
\end{eqnarray*}
Then we have the following proposition.
\begin{prop}\label{2pprop2}
The map $f\longmapsto \rho$ initially defined on $U(\mathscr{A})$ can be extended to an isometric isomorphism of $V(\mathscr{A})$ onto $\bigoplus\limits_{j\in\mathbb{N}}L^{2}((0,1],\ell^{2}
(\mathbb{Z}^{2n}),G_{0,0}^{\varphi_{j}})$.
\end{prop}
\begin{defn}
For $\{\varphi_{j}:j\in\mathbb{N}\}\subset\mathscr{A}$ and $\lambda\in(0,1], J(\lambda)$ is defined to be $\overline{span}\{\mathscr{T}\varphi_{j}(\lambda):j\in\mathbb{N}\}$.
\end{defn}
\begin{theorem}\label{2pth3}
Assume that each $\varphi_{j}$ satisfies condition C. Then $E(\mathscr{A})$ is a Parseval frame for $V(\mathscr{A})$ if and only if $\{\mathscr{T}\varphi_{j}(\lambda):j\in\mathbb{N}\}$ is a Parseval frame for $J(\lambda)$ for a.e. $\lambda\in\bigcap\limits_{j\in\mathbb{N}}\Omega_{\varphi_{j}}$, where $\Omega_{\varphi_{j}}=\{\lambda\in (0,1]:G_{0,0}^{\varphi_{j}}(\lambda)\neq 0\}$.
\end{theorem}
\begin{rem}\label{2prem3}
Fix $j\in\mathbb{N}$. Let $\varphi_{j}\in\mathscr{A}$ satisfy condition C. Then $\{L_{(2k,l,m)}\varphi_{j}:(k,l,m)\in\mathbb{Z}^{2n+1}\}$ is a Parseval frame for $V(\varphi_{j})$ if and only if $G_{0,0}^{\varphi_{j}}(\lambda)=1$ a.e. $\lambda\in\Omega_{\varphi_{j}}$.
\end{rem}
\begin{cor}\label{2pcor2}
Fix $j\in\mathbb{N}$. Let $\varphi_{j}\in\mathscr{A}$ satisfy condition C. Let $\psi_{j}\in L^{2}(\mathbb{H}^{n})$ be defined by
\begin{eqnarray*}
\widehat{\psi_{j}}(\lambda)&=\left\{
\begin{array}{l l}
G_{0,0}^{\varphi_{j}}(\lambda)^{-\frac{1}{2}}~\widehat{\varphi_{j}}(\lambda),& \quad \text{$\lambda\in\Omega_{\varphi_{j}}$,}\\
0, & \quad \text{otherwise.}
\end{array}\right.\label{2peq33}
\end{eqnarray*}
Then $\{L_{(2k,l,m)}\psi_{j}:(k,l,m)\in\mathbb{Z}^{2n+1}\}$ is a Parseval frame for $V(\varphi_{j})$.
\end{cor}
This result helps us to obtain a decomposition theorem for a shift-invariant space V in $ L^{2}(\mathbb{H}^{n})$.
\begin{defn}
If the collection $\{L_{(2k,l,m)}\varphi:(k,l,m)\in\mathbb{Z}^{2n+1}\}$ is a Parseval frame for $V(\varphi)$, then $\varphi$ is called Parseval frame generator of $V(\varphi)$.
\end{defn}
\begin{thm}
If V is a shift-invariant space in $L^{2}(\mathbb {H}^{n})$, then there exists a family of functions $\{\varphi_{\alpha}\}_{\alpha\in I}$ in $L^{2}(\mathbb {H}^{n})$ (where I is an index set) such that $V= \bigoplus\limits_{\alpha\in I} V(\varphi_{\alpha})$. In addition, if all the $\varphi_{\alpha}$ satisfy condition C, then there exists a family of functions $\{\psi_{\alpha}\}_{\alpha\in I}$ in $L^{2}(\mathbb {H}^{n})$ such that for each $\alpha\in I$, $\psi_{\alpha}$ is a Parseval frame generator of $V(\varphi_{\alpha})$. Moreover, in this case, if $f\in V$, then $\|f\|^{2}=\sum\limits_{\alpha\in I}\|\fontdimen16\textfont2=5pt\rho_{\alpha}\|^{2}_{L^{2}((0,1], \ell^{2}(\mathbb{Z}^{2n}),G^{\varphi_{\alpha}}_{0,0})}$, where $\fontdimen16\textfont2=5pt\rho_{\alpha}\in L^{2}((0,1], \ell^{2}(\mathbb{Z}^{2n}),G^{\varphi_{\alpha}}_{0,0})$.
\end{thm}
Similarly, we get a characterization for $E(\mathscr{A})$ to be a frame for $V(\mathscr{A})$. 
\begin{theorem}\label{2pth4}
Assume that each $\varphi_{j}$ satisfies condition C. Then $E(\mathscr{A})$ is a frame for $V(\mathscr{A})$ with frame bounds A,B if and only if $\{\mathscr{T}\varphi_{j}(\lambda):j\in\mathbb{N}\}$ is a frame for $J(\lambda)$ with frame bounds A, B for a.e. $\lambda\in\bigcap\limits_{j\in\mathbb{N}}\Omega_{\varphi_{j}}$.
\end{theorem}
\begin{theorem}
Assume that each $\varphi_{j}$ satisfies condition C. Then $E(\mathscr{A})$ is a Riesz basis for $V(\mathscr{A})$ with bound A, B if and only if $\{\mathscr{T}\varphi_{j}(\lambda):j\in\mathbb{N}\}$ is a Riesz basis for $J(\lambda)$ with bounds A, B for a.e. $\lambda\in (0,1]$.
\end{theorem}

\section{More on frames for $V(\varphi)$}
In this section, we show that if $\varphi\in L^{2}(\mathbb{H}^{n})$ has compact support, then the system consisting of left translates of $\varphi$ is a Bessel sequence and it is a non-redundant frame. Given two Bessel sequences $\{L_{(2k,l,m)}\varphi:(k,l,m)\in\mathbb{Z}^{2n+1}\}$ and $\{L_{(2k,l,m)}\tilde{\varphi}:(k,l,m)\in\mathbb{Z}^{2n+1}\}$, we provide a necessary and sufficient condition on the generators $\varphi$ and $\tilde{\varphi}$ such that $\{L_{(2k,l,m)}\tilde{\varphi}:(k,l,m)\in\mathbb{Z}^{2n+1}\}$ is an oblique dual of
$\{L_{(2k,l,m)}\varphi:(k,l,m)\in\mathbb{Z}^{2n+1}\}$ under a certain additional assumption.  If, in addition we assume that $\varphi,\tilde{\varphi}$ have compact support, we obtain much simpler characterization for the collection  $\{L_{(2k,l,m)}\tilde{\varphi}:(k,l,m)\in\mathbb{Z}^{2n+1}\}$ to be an oblique dual of
$\{L_{(2k,l,m)}\varphi:(k,l,m)\in\mathbb{Z}^{2n+1}\}$. We also show that if $\{L_{(2k,l,m)}\varphi:(k,l,m)\in\mathbb{Z}^{2n+1}\}$ is a frame for $V(\varphi)$, then there exists $\tilde{\varphi}\in L^{2}(\mathbb{H}^{n})$ such that $\{L_{(2k,l,m)}\tilde{\varphi}:(k,l,m)\in\mathbb{Z}^{2n+1}\}$ is an oblique dual of
$\{L_{(2k,l,m)}\varphi:(k,l,m)\in\mathbb{Z}^{2n+1}\}$. However, we are not able to show the uniqueness of $\tilde{\varphi}$.
\begin{prop}\label{2pprop3}
Let $\varphi\in L^{2}(\mathbb{H}^{n})$ have compact support satisfying condition C. Then the following holds:
\begin{itemize}
\item [(i)] $\{L_{(2k,l,m)}\varphi: (k,l,m)\in \mathbb{Z}^{2n+1}\}$ is a Bessel sequence.
\item [(ii)] $\{L_{(2k,l,m)}\varphi: (k,l,m)\in \mathbb{Z}^{2n+1}\}$ can not be an overcomplete frame sequence.
\end{itemize}
\end{prop}
\begin{theorem}\label{2pth5}
Let $\varphi, \tilde{\varphi}\in L^{2}(\mathbb{H}^{n})$ be such that they satisfy the following condition
\begin{eqnarray}\label{2peq25}
\sum\limits_{r\in\mathbb{Z}}\langle\widehat{\varphi}(\lambda+r),
\widehat{L_{(2k,l,0)}\tilde{\varphi}}(\lambda+r)\rangle_{\mathcal{B}_{2}}
|\lambda+r|^{n}=0~a.e.~\lambda\in (0,1],~\forall (k,l)\neq (0,0).
\end{eqnarray}
Assume that $\{L_{(2k,l,m)}\varphi: (k,l,m)\in \mathbb{Z}^{2n+1}\}$ and $\{L_{(2k,l,m)}\tilde{\varphi}: (k,l,m)\in \mathbb{Z}^{2n+1}\}$ are two Bessel sequences in $L^{2}(\mathbb{H}^{n})$. Then the following are equivalent:
\begin{itemize}
\item [(i)] $f=\sum\limits_{(k,l,m)\in\mathbb{Z}^{2n+1}}
\langle f,L_{(2k,l,m)}\tilde{\varphi}\rangle L_{(2k,l,m)}\varphi,~\forall~ f\in V(\varphi)$;
\item [(ii)] $\sum\limits_{r\in\mathbb{Z}}\langle \widehat{\varphi}(\lambda+r), \widehat{\tilde{\varphi}}(\lambda+r)\rangle_{\mathcal{B}_{2}}|\lambda+r|^{n}=1$ a.e. $\{\lambda: G_{0,0}^{\varphi}(\lambda)\neq 0\}$.
\end{itemize}
\end{theorem}
\begin{rem}\label{2pcor1}
If one of the equivalent conditions of Theorem \ref{2pth5} holds, then $\{L_{(2k,l,m)}\varphi\\: (k,l,m)\in \mathbb{Z}^{2n+1}\}$ is a frame for $V(\varphi)$.
\end{rem}
\begin{theorem}
Let $\varphi,\tilde{\varphi}\in L^{2}(\mathbb{H}^{n})$ be such that they satisfy  (\ref{2peq25}). Assume that $\varphi,\tilde{\varphi}$ have compact support and both of them satisfy condition C. Then the following are equivalent:
\begin{itemize}
\item [(i)] $f=\sum\limits_{(k,l,m)\in\mathbb{Z}^{2n+1}}
\langle f,L_{(2k,l,m)}\tilde{\varphi}\rangle L_{(2k,l,m)}\varphi,~\forall~ f\in V(\varphi)$;
\item [(ii)] $\langle \varphi, L_{(2k,l,m)}\tilde{\varphi}\rangle=\delta_{k,l,m,0}$.
\end{itemize}
\end{theorem}

\begin{theorem}
Let $\varphi\in L^{2}(\mathbb{H}^{n})$ and assume that $\{L_{(2k,l,m)}\varphi: (k,l,m)\in \mathbb{Z}^{2n+1}\}$ is a frame for $V(\varphi)$. Then there exists a function $\tilde{\varphi}\in V(\varphi)$ such that
\begin{eqnarray*}
f=\sum\limits_{(k,l,m)\in\mathbb{Z}^{2n+1}}
\langle f,L_{(2k,l,m)}\tilde{\varphi}\rangle L_{(2k,l,m)}\varphi, ~\forall f\in V(\varphi),
\end{eqnarray*}
namely $\tilde{\varphi}=S^{-1}\varphi$, where S is the frame operator corresponding to the frame $\{L_{(2k,l,m)}\varphi:(k,l,m)\in\mathbb{Z}^{2n+1}\}$.
\end{theorem}

\end{document}